\documentclass[12pt,a4paper]{article}

\usepackage{graphicx}
\usepackage{amsmath}
\usepackage{amsthm}
\usepackage{amssymb}
\usepackage{amsfonts}
\usepackage{latexsym}
\newcommand{\beql}[1]{\begin{equation}\label{#1}}
\newcommand{\eeq}{\end{equation}}
\newcommand{\comment}[1]{}

\newcommand{\Norm}[1]{{\left\|{#1}\right\|}}

\newcommand{\PP}{{\mathcal P}}
\newcommand{\PK}{{\mathcal P}_n(K)}

\newcommand{\sj}{{\frac{\sin\varphi_j}{r_j}}}

\newcommand{\Z}{{\mathcal Z}}
\newcommand{\W}{{\mathcal W}}
\newcommand{\RR}{{\mathbb R}}
\newcommand{\CC}{{\mathbb C}}

\newcommand{\NN}{{\mathbb N}}

\newcommand{\dist}{{\rm dist\,}}
\newcommand{\con}{{\rm con\,}}
\newcommand{\diam}{{\rm diam\,}}
\newcommand{\width}{{\rm width\,}}

\newcommand{\bnu}{\mbox{\boldmath$\nu$}}

\newcounter{dfn}
\def\thedfn{\arabic{dfn}}

\newcounter{thm}
\newcounter{othm}
\def\theothm{\Alph{othm}}
\newenvironment{othm}{
  \em
  \vskip 0.10in
  \refstepcounter{othm}
  \noindent{\bf Theorem\ \theothm}
}{\vskip 0.10in}

\newenvironment{olemma}{
  \em
  \vskip 0.10in
  \refstepcounter{othm}
  \noindent{\bf Lemma\ \theothm}
}{\vskip 0.10in}

\newcounter{mysec}
\def\themysec{\arabic{mysec}}
\newcommand{\mysection}[1]{
  \vskip 0.25in
  \refstepcounter{mysec}\centerline{\large\bf \S\themysec.\ {#1}}\par
  \nopagebreak
  \addcontentsline{toc}{section}{{\bf \themysec.}\ {#1}}
}
\newcounter{mysubsec}[mysec]
\newtheorem{theorem}{Theorem}

\newtheorem{lemma}{Lemma}

\theoremstyle{definition}

\newtheorem{remark}{Remark}

\newcommand{\de}{\delta}

\reversemarginpar 

\begin{document}

\title{Right order Tur\'an-type converse
Markov inequalities for convex domains on the plane}

\author{Szil\'ard Gy. R\'ev\'esz}

\maketitle

\begin{abstract}
For a convex domain $K\subset\CC$ the well-known general
Bernstein-Markov inequality holds asserting that a polynomial $p$
of degree $n$ must have $\Norm{p'} \le c(K) n^2 \Norm{p}$. On the
other hand for polynomials in general, $\Norm{p'}$ can be
arbitrarily small as compared to $\Norm{p}$.

The situation changes when we assume that the polynomials in
question have all their zeroes in the convex body $K$. This was
first investigated by Tur\'an, who showed the lower bounds
$\Norm{p'} \ge (n/2) \Norm{p}$ for the unit disk $D$ and
$\Norm{p'} \ge c \sqrt{n} \Norm{p}$ for the unit interval
$I:=[-1,1]$. Although partial results provided general lower
estimates of lower order, as well as certain classes of domains
with lower bounds of order $n$, it was not clear what order of
magnitude the general convex domains may admit here.

Here we show that for all compact and convex domains $K$ with
nonempty interior and polynomials $p$ with all their zeroes in
$K$ $\Norm{p'} \ge c(K) n \Norm{p}$ holds true, while $\Norm{p'}
\le C(K) n \Norm{p}$ occurs for any $K$. Actually, we determine
$c(K)$ and $C(K)$ within a factor of absolute numerical
constant.
\end{abstract}
\let\oldfootnote\thefootnote
\def\thefootnote{}
\footnotetext{Supported in part in the framework of the
Hungarian-French Scientific and Technological Governmental
Cooperation, Project \# F-10/04 and the Hungarian-Spanish
Scientific and Technological Governmental Cooperation, Project \#
E-38/04.}
\footnotetext{This work was not supported by Hungarian
National Foundation for Scientific Research.}
\let\thefootnote\oldfootnote

\bigskip
\bigskip
\bigskip
\bigskip
\bigskip
\bigskip
\bigskip
\bigskip
\bigskip



\newpage

%%%%%%%%%%%%%%%%%%%%%%%%%%%%%%%%%%%%%%%%%%%%%
%%%%%%%%%%         Introduction
%%%%%%%%%%%%%%%%%%%%%%%%%%%%%%%%%%%%%%%%%%%%%

\mysection{Introduction} \label{sec:introduction}

\bigskip

%%%%%%%%%%%%%%%%%%%%%%%%%%%%%%%%%%%%%%%%%%%%%
%%%%%%%%%        Initial work of Turan
%%%%%%%%%%%%%%%%%%%%%%%%%%%%%%%%%%%%%%%%%%%%%

\noindent On the complex plane polynomials of degree $n$ admit a
Bernstein-Markov inequality\footnote{Namely, to each point $z$ of
$K$ there exists another $w\in K$ with $|w-z|\ge {\rm diam}(K)/2$,
and applying Markov's inequality on the segment $[z,w]\subset K$
yields $|p'(z)| \le  (1/{\rm diam}(K)) n^2 \Norm{p}_K$.}
$\Norm{p'}_K\le c_K n^2 \Norm{p}_K$ on all convex, compact
$K\subset \CC$. Here the norm $\Norm{\cdot} :=\Norm{\cdot}_K$
denotes sup norm over values attained on $K$.

Sixty-five years ago Paul Tur\'an studied converse inequalities
of the form $\Norm{p'}_K\ge c_K n^A \Norm{p}_K$. Clearly such a
converse can hold only if further restrictions are imposed
on the occurring polynomials $p$. Tur\'an assumed that all zeroes
of the polynomials must belong to $K$. So denote the set of
complex (algebraic) polynomials of degree (exactly) $n$ as
$\PP_n$, and the subset with all the $n$ (complex) roots in some
set $K\subset\CC$ by $\PK$. The (normalized) quantity under our
study is thus the "inverse Markov factor"
\begin{equation}\label{Mdef}
M_n(K):=\inf_{p\in \PK} M(p) \qquad \text{\rm with} \qquad
M:=M(p):=\frac{\Norm{p'}}{\Norm{p}}~.
\end{equation}

\begin{othm}{\bf[Tur\'an].}\label{oth:Turandisk} If $p\in \PP_n(D)$, where
$D$ is the unit disk, then we have
\begin{equation}\label{Turandisk}
\Norm{p'}_D\ge \frac n2 \Norm{p}_D~.
\end{equation}
\end{othm}

\begin{othm}{\bf[Tur\'an].}\label{oth:Turanint} If $p\in\PP_n(I)$,
where $I:=[-1,1]$, then we have
\begin{equation}\label{Turanint}
\Norm{p'}_I\ge \frac {\sqrt{n}}{6} \Norm{p}_I~.
\end{equation}
\end{othm}

Theorem \ref{oth:Turandisk} is best possible. Regarding Theorem
\ref{oth:Turanint}, Tur\'an pointed out by example of
$(1-x^2)^{n}$ that the $\sqrt{n}$ order is sharp. The slightly
improved constant $1/(2e)$ can be found in \cite{LP}, and the
value of the constant is computed for all fixed $n$ precisely in
\cite{Er}.

The key to Theorem \ref{oth:Turandisk} was the following
observation, which had already been present implicitly in
\cite{Tur} and \cite{Er} and was later formulated explicitly in
\cite{LP}.

\begin{olemma}
{\bf[Tur\'an].}\label{Tlemma} Assume that $z\in\partial K$ and
that there exists a disc $D_R$ of radius $R$ so that
$z\in\partial D_R$ and $K\subset D_R$. Then for all $p\in\PK$ we
have
\begin{equation}\label{Rdisc}
|p'(z)| \ge \frac n{2R} |p(z)|~.
\end{equation}
\end{olemma}

Drawing from the work of Tur\'an, Er\H od \cite{Er} already
addressed the question: "For what kind of domains does the method
of Tur\'an provide $cn$ order of oscillation for the derivative?"
In particular, he showed

\begin{othm}{\bf[Er\H od].}\label{th:ellipse}
Let $0<b<1$ and let $E_b$ denote the ellipse domain with major
axes $[-1,1]$ and minor axes $[-ib,ib]$. Then for all
$p\in\PP_n(E_b)$ we have
\begin{equation}\label{oldellipse}
\Norm{p'} \ge \frac {b}{2}  n \Norm{p}.
\end{equation}
\end{othm}

Moreover, he elaborated on the inverse Markov factors belonging to
domains with some favorable geometric properties, such as having
positive curvature exceeding a given fixed positive bound at all
boundary points, or at all boundary points with the exception of
a given (finite) set of vertices, etc. For a detailed account of
results of Er\H od in this direction, as well as even further
results applying basically Tur\'an's Lemma \ref{Tlemma}, see the
recent works \cite{LP}, \cite{E} and \cite{R}.

A lower estimate of the inverse Markov factor for any convex set
and of at least the same order as for the interval was obtained
in full generality only in about three years ago.

\begin{othm}{\bf[Levenberg-Poletsky].}\label{th:generalroot} If
$K\subset \CC$ is a compact, convex set, $d:=\diam{K}$ and $p\in
\PK$, then we have
\begin{equation}\label{genrootineq}
\Norm{p'}\ge \frac {\sqrt{n}}{20\,\diam(K)}  \Norm{p}~.
\end{equation}
\end{othm}

Interestingly, it turned out that among all convex compacta only
intervals can have an inverse Markov constant of such a small
order, while for convex compact domains $K$ and for all $p\in\PK$
we have at least $M(p)\ge C_1(K) n^{2/3}$, see \cite{Rev}. Recall
that here the term {\em convex domain} stands for a compact,
convex subset of $\CC$ {\em having nonempty interior}. Clearly,
assuming boundedness is natural, since all polynomials of
positive degree have $\Norm{p}_K=\infty$ when the set $K$ is
unbounded. Also, all convex sets with nonempty interior are {\em
fat}, meaning that ${\rm cl}(K)={\rm cl}({\rm int} K)$. Hence
taking the closure does not change the $\sup $ norm of
polynomials under study. The only convex, compact sets, falling
out by our restrictions, are the intervals, for what Tur\'an has
already shown that his lower estimate is of the right order.

The case of the unit disk and the example of $p(z)=1+z^n$ shows
that in general the order of the inverse Markov factor can not be
higher than $n$. On the other hand, some general classes of
domains were found to have order $n$ inverse Markov factors. Let
us list a few examples of such domains.

\begin{enumerate}
  \item All convex domains with $C^2$-smooth boundary and
  curvature above a given fixed parameter $\kappa>0$ (Er\H od
  \cite{Er}\footnote{For further details on the constants and
  more precise details of the slightly incomplete proofs in
  \cite{Er} for items 1-3 see \cite{R}.} and R\'ev\'esz \cite{R}).
  \item Convex domains bounded by finitely many $C^2$-smooth
  Jordan arcs and a finite number of vertices, with the curvature
  of any relative interior points of the arcs bounded away from 0
  (Er\H od \cite{Er} and R\'ev\'esz \cite{R}).
  \item Convex domains of smooth boundary and curvature bounded
  away from 0, with the exception of one straight line segment on
  the boundary having length $< \diam (K)/4$, (Er\H od \cite{Er}).
  \item A square\footnote{Erd\'elyi also proves similar results on
  rhombuses, under the further condition of some symmetry of the
  polynomials in consideration -- e.g. if the polynomials are real,
  or odd. Note also that his work \cite{E} preceded \cite{Rev, R}
  and apparently was accomplished without being aware of details of
  \cite{Er}.} (Erd\'elyi, \cite{E}, \cite{E2}).
  \item Convex domains with finitely many vertices having vertices
  of only acute supplementary angles and finitely many smooth
  Jordan arcs connecting the vertices \cite{Rev, R}.
  \item Smooth convex domains \cite{Rev, R}.
  \item Convex domains of {\em fixed positive depth} \cite{Rev, R}.
  \item Convex domains with their almost everywhere (with respect
  to arc length measure) existing curvature exceeding almost everywhere
  a given positive lower bound \cite{R}.
\end{enumerate}

For further details and a discussion of the results of Er\H od
\cite{Er}, see the references, in particular \cite{R}. On the
other hand, it was not known whether the inverse Markov factor
can be $o(n)$ or not.

To study \eqref{Mdef} some geometric parameters of the convex
domain $K$ are involved naturally. We write $d:=d(K):=\diam (K)$
for the {\em diameter} of $K$, and $w:=w(K):={\width}(K)$ for the
{\em minimal width} of $K$. Note that a convex domain is a
closed, bounded, convex set $K\subset\CC$ with nonempty interior,
hence $0<w(K)\le d(K)<\infty$. Our main result is the following.

\begin{theorem}\label{th:main} Let $K\subset\CC$ be any
convex domain. Then for all $p\in\PK$ we have
\begin{equation}\label{mainresult}
\frac{\Norm{p'}}{\Norm{p}} \ge C(K) n \qquad \text{\rm with}
\qquad C(K)=0.0006 \frac{w(K)}{d^2(K)}~.
\end{equation}
\end{theorem}

Clearly this result contains all the above results apart from the
precise value of the absolute constant factor. Moreover, the
result is essentially sharp for all convex domains $K$: see \S
\ref{sec:sharpness} below.

%%%%%%%%%%%%%%%%%%%%%%%%%%%%%%%%%%%%%%%%%%%%%%%%%%%%%
%%%%%%%%     The Proof
%%%%%%%%%%%%%%%%%%%%%%%%%%%%%%%%%%%%%%%%%%%%%%%%%%%%%

\mysection{Proof of Theorem \ref{th:main}}\label{sec:proof}
\bigskip
Our proof will follow the argument of \cite{Rev}, with one key
alteration, suggested to us by G\'abor Hal\'asz. Namely, we start
with picking up a boundary point $\zeta\in\partial K$ of
maximality of $|p|$, and consider a supporting line at $\zeta$ to
$K$ as in \cite{Rev}. However, then we do not use the {\em
normal} direction to compare values of $p$ at $\zeta$ and on the
intersection of $K$ and this normal line, but instead here we
compare the values of $p$ at $\zeta$ and {\em on a line slightly
slanted off from the normal}. Comparing the calculations here and
in \cite{Rev} the reader will detect how this change led to an
essential improvement of the result through improving the
contribution of the factors belonging to zeroes close to the
supporting line. In \cite{Rev} we could get a square term (in $h$
there) only, due to orthogonality and the consequent use of the
Pithagorean Theorem in calculating the distances. However, here
we obtain {\em linear dependence} in $\delta$ via the general
cosine theorem for the slanted segment $J$. (That insightful
observation was provided by G. Hal\'asz.) One of the major
geometric features still at our help is the fact, that {\em when
the intersection} of a normal or close-to-normal line with $K$
{\em is small}, {\em then one part of the convex domain} $K$, cut
into two by the line, {\em will also be small in the same order}.
That was explicitly formulated in \cite{Rev}, and is used
implicitly even here through various calculations with the
angles: this is the key feature which allows us to bend the
direction of the normal a bit {\em towards the smaller portion}
of $K$. As a result of the improved estimates squeezed out this
way, we do not need to employ the second technique, also going
back to Tur\'an, i.e. integration of $(p'/p)'$ over a suitably
chosen interval. As pointed out already in \cite{Rev}, this part
of the proof yields weaker estimates than $cn$, so avoiding it is
not only a matter of convenience, but is an essential necessity.

\begin{proof} We list the zeroes of a polynomial $p\in\PK$ according to
multiplicities as $z_1,\dots,z_n$, and the set of these zero
points is denoted as $\Z:=\Z(p):=\{ z_j~:~ j=1,\dots,n\}\subset
K$. (It suffices to assume that all $z_j$ are distinct, so we do
not bother with repeatedly explaining multiplicities, etc.)
Assume, as we may, $p(z)=\prod_{j=1}^n (z-z_j)$.

We start with picking up a point $\zeta$ of $K$, where $p$
attains its norm. By the maximum principle, $\zeta\in\partial K$,
and by convexity there exists a supporting line to $K$ at $\zeta$
with inward normal vector $\bnu$, say. Without loss of generality
we can take $\zeta=0$ and $\bnu=i$. Now by definition of
the minimal width $w=w(K)$, there exists a point $A\in K$ with
$\Im A \ge w$; by symmetry, we may assume $\Re A \le 0$, say.

Sometimes we write the zeroes in their polar form
\begin{equation}\label{polarform}
z_j=r_je^{i\varphi_j} \qquad \left( r_j:=|z_j|, \quad
\varphi_j:=\arg z_j \quad (j=1,\dots,n) \right)~.
\end{equation}
Throughout the proofs with $[(\varphi,\psi)]$ being any open,
closed, halfopen-half-closed or halfclosed-halfopen interval we
use the notations
\begin{equation}\label{Ssector}
S[(\varphi,\psi)]:=\{z\in\CC~:~ \arg(z)\in [(\varphi,\psi)]\}
\end{equation}
and
\begin{equation}\label{Znsector}
\Z[(\varphi,\psi)]:=\Z\cap S[(\varphi,\psi)]\,, \qquad
n[(\varphi,\psi)]:=\# \Z[(\varphi,\psi)]\,,
\end{equation}
for the sectors, the zeroes in the sectors, and the number of
zeroes in the sectors determined by the angles $\varphi$ and
$\psi$.

Let us formulate a well-known but useful fact in advance.

\begin{lemma}\label{l:capacity} Let $J=[u,v]$ be any interval on
the complex plane with $u\ne v$ and let $J \subset R \subset \CC$
be any set containing $J$. Then for all $k\in\NN$ we have
\begin{equation}\label{capacity}
\min_{w_1,\dots,w_k\in R} \max_{z\in J} \left| \prod_{j=1}^k
(z-w_j) \right| \ge 2 \left(\frac{|J|}{4}\right)^k~.
\end{equation}
\end{lemma}

\begin{proof} This is essentially the classical result of Chebyshev
for a real interval, cf. \cite{BE, MMR}, and it holds for much
more general situations (perhaps with the loss of the factor 2)
from the notion of Chebyshev constants and capacity, cf. Theorem
5.5.4. (a) in \cite{Rans}.
\end{proof}

In all our proof we fix the angles
\begin{equation}\label{fidef}
\psi:=\arctan \left( \frac{w}{d} \right)\in (0,\pi/4] ~\qquad
\text{\rm and} \qquad \theta:=\psi/20 \in (0,\pi/80].
\end{equation}

Since $|p(0)|=\|p\|$, $M\ge |p'(0)/p(0)|$. Observe that for any
subset $\W\subset\Z$ we then have
\begin{equation}\label{Impartsum}
M\ge \left| \frac{p'}{p}(0) \right| \ge -\Im {\frac{p'}{p}(0)} =
\sum_{j=1}^n \Im \frac {-1}{z_j} \ge \sum_{z_j\in\W} \Im \frac
{-1}{z_j}=\sum_{z_j\in\W} \frac {\sin \varphi_j}{r_j} \,,
\end{equation}
since all terms in the full sum are nonnegative.

Let us consider now the ray (straight half-line) emanating from
$\zeta=0$ in the direction of $e^{i(\pi/2-2\theta)}$. This ray
intersects $K$ in a line segment $[\zeta,D]$, and if $D=\zeta$,
then $K\subset S[\pi/2-2\theta,\pi]$ and a standard argument using
e.g. Tur\'an's Lemma \ref{Tlemma} yields $M\ge n/(2d)$. Hence we
may assume $D\ne 0$.

Consider now any point $B\in K$ with maximal real part, and take
$B':=\Re B=\max \{ \Re z ~:~ z\in K\}$. Since $D \ne 0$, $B'> 0$,
and as $\Re A\le 0$ and $\Re B$ is maximal, $[A,B']$ intersects
$[0,D]$ in a point $D'\in [0,D]$, i.e. $[0,D']\subset
[0,D]\subset K$. Moreover, the angle at $B'$ between the real line
and $AB'$ is $-\arg(B'-A)=-\arg(B'-D') \in [\psi, \pi/2)$. Indeed,
$\Im (A-B')\ge w$ and $\Re (B'-A)=\Re(B-A)\le d$ (resulting from
$A,B\in K$) imply $-\arg(B'-A)\ge \arctan (w/d)=\psi$.

In the following let us write $\delta:=|D'|>0$; it can not
vanish, as $B'\ne 0$ and the line segment $[B',A]$ intersects the
real line only in $B'$. Consider the point $B"\in\RR$ with $B"\ge
B' >0$ and $-\arg(B"-D')=\psi$. We can say now that K lies both
in the upper half of the disk with radius $d$ around 0 (which we
denote by $U$), and the halfplane $\Re z \le B"$ (which we denote
by $H$); moreover, $[0,D']\subset K \subset (U\cap H)$.

Now we put $D":=3D'/4$ and take
\begin{equation}\label{Jdef}
J:= \left[D", D' \right] \subset K \quad \text{i.e.}\;\;
J:=\{\tau:=te^{i(\pi/2-2\theta)}\de~:~3/4\le t \le 1\}~.
\end{equation}

Denoting $D_r(0):=\{z\,:\, |z|\le r\}$ we split the set $\Z$ into
the following parts.
  \begin{align}\label{Zsplitup}
\Z_1:&= \Z[0,\theta]\,, \qquad \qquad \qquad \qquad
\qquad\qquad\qquad \qquad \quad \mu:=\#\Z_1=n[0,\theta] \notag \\ %%%
\Z_{2}:&= \Z(\theta,\pi-\theta)\cap \left\{\Im (e^{i2\theta}z)<
\frac 38 \delta \right\}\,,\qquad \qquad
\qquad\qquad\qquad \nu:=\#\Z_{2} \notag \\ %%%
\Z_{3}:&= \Z(\theta,\pi-\theta)\cap\left\{\Im (e^{i2\theta}z) \ge
\frac 38 \delta \right\} \cap D_{2\delta}(0)\,,\qquad \qquad
\qquad~ \kappa:=\#\Z_{3} \notag \\ %%%
\Z_{4}:&= \Z(\theta,\pi-\theta)\cap\left\{\Im (e^{i2\theta}z) \ge
\frac 38 \delta \right\}\setminus D_{2\delta}(0)= \\ %%%
    &=\Z(\theta,\pi- \theta)\setminus \Z_{2}\setminus \Z_{3}\,,
    \qquad\qquad\qquad\qquad\qquad\qquad\qquad\quad k:=\#\Z_{4} \notag \\ %%%
\Z_5:&= \Z[\pi-\theta,\pi]\,,\qquad\qquad\qquad\qquad \qquad\qquad
\quad m:=\#\Z_5=n[\pi-\theta,\pi] \notag ~.
  \end{align}

In the following we establish an inequality from condition of
maximality of $|p(0)|$. First we estimate the distance of any
$z_j\in\Z_1$ from $J$. In fact, taking any point
$z=re^{i\varphi}\in H\cap S[0,\theta]$ the sine theorem yields
$r\cos \varphi = \Re{z} \le |B"|=\delta(\sin(\pi/2+2\theta-\psi)/
\sin\psi) = \delta\cos(\psi-2\theta)/\sin(\psi)< \de \cot(18
\theta)$, and so
\begin{equation}\label{deltarst}
r\sin\theta < \frac{\sin \theta}{\cos \varphi} \frac {\delta}{\tan
(18\theta)} \le \delta \frac{\tan\theta}{\tan(18\theta)} <
\frac{\delta}{18} ~.
\end{equation}

Now $\dist (z,J)=\min_{3/4\le t\le 1} |z-\tau|$, and by the
cosine theorem
$|z-\tau|^2=t^2\delta^2+r^2-2\cos(\pi/2-\varphi-2\theta)~rt\delta$.
Since
$\cos(\pi/2-\varphi-2\theta)=\sin(\varphi+2\theta)\le\sin(3\theta)\le
3 \sin\theta$, \eqref{deltarst} implies $|z-\tau|^2 \ge
t^2\delta^2+r^2 - 6t\delta \sin\theta ~r \ge t^2\delta^2+r^2 -
(1/3)t\delta^2 $, and thus $\min_{3/4\le t\le 1} |z-\tau|^2 \ge
\min_{3/4\le t\le 1} t^2\delta^2+r^2 - (1/3) t\delta^2 =
r^2+(5/16) \delta^2 $. It follows that we have
$$ %%%
\frac{|z-\tau|^2}{|z|^2}\ge \frac{r^2+(5/16)\delta^2}{r^2}
> 1+ \frac{(90/16)\sin\theta\,\,\delta}{r} > 1+
\frac{5\sin\theta\,\,\delta}{d} \quad \left(\tau\in J\right)\,,
$$ %%%
applying also \eqref{deltarst} to estimate $\de/r$ in the last but
one step. Now $\de/d\le 1$ and $5\sin\theta < 0.2$, hence we can
apply $\log(1+x)\ge x-x^2/2 \ge 0.9 x$ for $0<x<0.2$ to get
$$ %%%
\frac{|z-\tau|^2}{|z|^2} \ge \exp\left(0.9 \frac{5\sin\theta ~
\delta}{d}\right) > \exp\left(\frac{4 \sin\theta ~ \delta}{d}
\right) \quad \left(\tau \in J\right)~.
$$ %%%
Applying this estimate for all the $\mu$ zeroes $z_j\in\Z_1$ we
finally find
\begin{equation}\label{Z1zeros}
\prod_{z_j\in\Z_1}\left|\frac{z_j-\tau}{z_j}\right| \ge
\exp\left(\frac{2\sin\theta ~ \delta\mu}{d} \right)~\qquad
\left(\tau=t\de e^{i(\pi/2-2\theta)} \in J \right).
\end{equation}
The estimate of the contribution of zeroes from $\Z_5$ is somewhat
easier, as now the angle between $z_j$ and $\tau$ exceeds $\pi/2$.
By the cosine theorem again, we obtain for any $z=r
e^{i\varphi}\in S[\pi-\theta,\pi]\cap U$ the estimate
\begin{align}\label{franc}
|z-\tau|^2 = &
r^2+t^2\de^2-2\cos(\varphi-(\pi/2-2\theta))\,rt\delta \notag\\ \ge
& r^2+t^2\de^2+2\sin\theta ~ rt\delta >
r^2\left(1+\frac{3\sin\theta ~ \delta}{2d}\right) ~~~\left(\tau
\in J\right)~,
\end{align}
as $t\ge 3/4$ and $r\le d$. Hence using again $\de/d\le 1$ and
$1.5\sin\theta < 0.06$ we can apply $\log(1+x)\ge x-x^2/2 \ge 0.97
x$ for $0<x<0.06$ to get $$ \frac{|z-\tau|}{|z|} \ge
\exp\left(\frac12 0.97 \frac{3\sin\theta ~ \delta}{2d}\right) \geq
\exp\left(\frac{18\sin\theta\,\delta}{25d} \right) ~~~\left(\tau
\in J\right)~, $$ whence
\begin{equation}\label{Z4contri}
\prod_{z_j\in\Z_5}\left|\frac{z_j-\tau}{z_j}\right| \ge\exp\left(
\frac{18\sin\theta\,\,\delta m}{25d} \right) ~\qquad
\left(\tau=t\de e^{i(\pi/2-2\theta)} \in J \right).
\end{equation}
Observe that zeroes belonging to $\Z_{2}$ have the property that
they fall to the opposite side of the line
$\Im(e^{i2\theta}z)=3\de/8$ than $J$, hence they are closer to $0$
than to any point of $J$. It follows that
\begin{equation}\label{Zstarcontri}
\prod_{z_j\in\Z_{2}}\left|\frac{z_j-\tau}{z_j}\right| \ge 1
~\qquad \left(\tau=t\de e^{i(\pi/2-2\theta)} \in J \right).
\end{equation}
Next we use Lemma \ref{l:capacity} to estimate the contribution of
zero factors belonging to $\Z_{3}$. We find
\begin{equation}\label{Z2pcontri}
\max_{\tau\in J} \prod_{z_j\in\Z_{3}}
\left|\frac{z_j-\tau}{z_j}\right| \ge
2\left(\frac{|J|}{4}\right)^{\kappa^{}} \prod_{z_j\in\Z_{3}}
\frac{1}{r_j} >\left(\frac{1}{32}\right)^{\kappa^{}} >
\exp(-3.5\kappa^{})~,
\end{equation}
in view of $|J|=\de/4$ and $r_j \le 2 \de$.

Note that for any point $z=re^{i\varphi}\in D_{2\delta}(0)\cap
\{\Im(e^{i2\theta}z)\ge 3\de/8 \}$ we must have $$
\frac{3\delta}{8}\le
\Im(e^{i2\theta}re^{i\varphi})=r\sin(\varphi+2\theta)~, $$ hence
by $r\le 2\de$ also $$ \sin(\varphi+2\theta)\ge \frac{3\delta}{8r}
\ge \frac 3{16} $$ and $\sin \varphi\ge
\sin(\varphi+2\theta)-2\theta \ge 3/16 - \pi/40 > 1/10$. Applying
this for all the zeroes $z_j\in \Z_{3}$ we are led to
\begin{equation}\label{maxcontri}
1\le \frac{2\de}{r_j} \le 20 \de \sj \qquad
\left(z_j\in\Z_{3}\right)~.
\end{equation}
On combining \eqref{Z2pcontri} with \eqref{maxcontri} we are led
to
\begin{equation}\label{Z2plus}
\max_{\tau\in J} \prod_{z_j\in\Z_{3}}
\left|\frac{z_j-\tau}{z_j}\right| \ge \exp\left(-70\de
\sum_{z_j\in\Z_{3}}\sj \right)~.
\end{equation}

Finally we consider the contribution of the zeroes from $\Z_{4}$,
i.e. the "far" zeroes for which we have
$\Im(z_je^{2i\theta})\ge3\de/8$, $\varphi_j\in
(\theta,\pi-\theta)$ and $|r_j|\ge 2\de$. Put now
$w:=z_je^{2i\theta}=u+iv=re^{i\alpha}$, and $s:=|\tau|=t\de$, say.
We then have
\begin{align}\label{farupzeroes}
\left|\frac{z_j-\tau}{z_j}\right|^2 &= \frac{|w-t\de i|^2}{r^2}=
\frac{u^2+(v-s)^2}{r^2}=1-\frac{2vs}{r^2}+\frac{s^2}{r^2} \\ &>
1-\frac{2vs}{r^2}+\frac{s^2}{r^2}\frac{v^2}{r^2} =
\left(1-\frac{vs}{r^2}\right)^2\ge
\left(1-\frac{|v|\delta}{r^2}\right)^2=
\left(1-\frac{\de|\sin\alpha|}{r}\right)^2. \notag
\end{align}
Recall that $\log(1-x) >-x-\frac{x^2}{2}\frac {1}{1-x}\ge
-x(1+1/2)$ whenever $0\le x \le 1/2$. We can apply this for
$x:=\de|\sin\alpha|/r_j\le \de/r_j \le 1/2$ using $r=r_j=|z_j|=|w|
\ge 2\de$. As a result, \eqref{farupzeroes} leads to
\begin{equation}\label{farupzeroescontr}
\left|\frac{z_j-\tau}{z_j}\right|\ge \exp\left( -\frac 32
\de\frac{|\sin(\varphi_j+2\theta)|}{r_j}\right)~,
\end{equation}
and using $|\sin(\varphi_j+2\theta)|\le
\sin(\varphi_j)+\sin(2\theta) \le 3 \sin\varphi_j$ (in view of
$\varphi_j\in(\theta,\pi-\theta)$), finally we get
\begin{equation}\label{Z3plusfin}
\prod_{z_j\in\Z_{4}}\left|\frac{z_j-\tau}{z_j}\right| \ge
 \exp\left( -\frac {9\de}{2} \sum_{z_j\in\Z_{4}}
 \frac{\sin\varphi_j}{r_j}\right)\qquad \left(\tau=t\de
e^{i(\pi/2-2\theta)} \in J \right)~.
\end{equation}

Collecting the estimates  \eqref{Z1zeros} \eqref{Z4contri}
\eqref{Zstarcontri} \eqref{Z2plus}  and \eqref{Z3plusfin} gives
for a certain point of maxima $\tau_0\in J$ in \eqref{Z2plus} the
inequality
\begin{align}\label{}
1 \ge & \frac{|p(\tau_0)|}{|p(0)|}= \prod_{z_j\in\Z
}\left|\frac{z_j-\tau_0}{z_j}\right| > \\ &
\exp\left\{\frac{18}{25} \sin\theta\,\delta\frac{\mu+m}{d}- 70\de
\sum_{z_j\in \Z_{2}\cup\Z_{3}\cup\Z_{4}} \sj \right\}~,\notag
\end{align}
or, after taking logarithms and cancelling by $18\de/25$
\begin{equation}\label{pointofmaxima}
\sin\theta \frac{\mu+m}{d} < \frac{875}{9} \sum_{z_j\in
\Z_{2}\cup\Z_{3}\cup\Z_{4}} \sj
\end{equation}
Observe that for the zeroes in $\Z_2\cup\Z_3\cup\Z_4$ we have
$\sin\varphi_j > \sin \theta$, whence also
\begin{equation}\label{triviesti}
(\nu+\kappa+k)\frac{\sin\theta}{d} \le \sum _{z_j\in
\Z_2\cup\Z_3\cup\Z_4} \sj~.
\end{equation}
Adding \eqref{pointofmaxima} and \eqref{triviesti} and taking into
account $\#\Z=\sum_{j=1}^5 \#Z_j$, we obtain
\begin{equation}\label{almostfinal}
\sin\theta \frac{n}{d}=\sin\theta \frac{\mu+m+\nu+\kappa+k}{d}<
\frac{884}{9} \sum _{z_j\in \Z_2\cup\Z_3\cup\Z_4} \sj~.
\end{equation}
Making use of \eqref{Impartsum} with the choice of $\W:=
\Z_2\cup\Z_3\cup\Z_4$ we arrive at
$$ %%%
\sin\theta \frac{n}{d}< \frac{884}{9} M ~,
$$ %%%
that is,
\begin{equation}\label{Mnesti}
M > \frac{9\sin\theta}{884d} n ~.
\end{equation}
It remains to recall \eqref{fidef} and to estimate
$$
\sin\theta = \sin \left( \frac{\arctan(w/d)}{20}\right)~.
$$
As $\theta\in (0,\pi/80]$, \; $\sin\theta >
\theta(1-\theta^2/6)\ge \theta(1-\pi/240)>0.98\theta$ and as
$0<w/d\le 1$, \; $\arctan(w/d)\ge (w/d)/(\pi/4)$, whence
$$ %%%
\sin\theta \ge 0.98 \frac{\arctan(w/d)}{20} \ge \frac{0.98}{5\pi}
\frac wd~.
$$ %%%
Substituting this last estimate into \eqref{Mnesti} yields
$$
M > \frac{9}{884} \cdot \frac{0.98}{5\pi}\cdot \frac w{d^2} \cdot
n
> 0.0006 \frac{w}{d^2} n~,
$$ concluding the proof. \end{proof}

\mysection{On sharpness of the main result}\label{sec:sharpness}

\begin{theorem}\label{th:sharp} Let $K\subset\CC$ be any
compact, connected set with diameter $d$ and minimal width $w$.
Then for all $n>n_0:=n_0(K):= 2 (d/16w)^2 \log (d/16w)$ there
exists a polynomial $p\in\PK$ of degree exactly $n$ satisfying
\begin{equation}\label{sharpresult}
\Norm{p'} \leq ~ C'(K)~ n~ {\Norm{p}} \qquad \text{\rm with}
\qquad C'(K):= 600 ~\frac{w(K)}{d^2(K)}~.
\end{equation}
\end{theorem}

\begin{remark} Note that here we do not assume that $K$ be convex,
but only that it is a connected, closed (compact) subset of
$\CC$. (Clearly the condition of boundedness is not restrictive,
$\|p\|$ being infinite otherwise.)
\end{remark}

\begin{proof} Take $a,b\in K$ with $|a-b|=d$ and $m\in\NN$ with
$m>m_0$ to be determined later. Consider the polynomials
$q(z):=(z-a)(z-b)$, $p(z)=(z-a)^m(z-b)^m=q^m(z)$ and
$P(z)=(z-a)^m(z-b)^{m+1}=(z-b)q^m(z)$. Clearly, $p,P\in \PK$ and
$\deg p= 2m$, $\deg P=2m+1$. We claim that these polynomials
satisfy inequality \eqref{sharpresult} for appropriate choice of
$m_0$.

First we make a few general observations. One obvious fact is
that if the unit vector $e:=(b-a)/d$, then the line
$\ell:=\{\frac{(a+b)}{2} + ite~:~t\in \RR\}$ separates $a$ and
$b$. Since $K$ is connected, also $\ell$ contains some point $c$
of $K$. Therefore, $\|q\|\geq |q(c)|=(d/2)^2+t^2 \geq (d/2)^2$.
Also, it is clear that $q'(z)=2z-a-b$ and hence $\|q'\| \leq
|z-a|+|z-b|\leq 2d$, by definition of the diameter.

As for $p$, we have $p'=mq'q^{m-1}$, hence
\begin{equation}\label{pprime}
\|p'\| \leq m \|q'\|  \|q\|^{m-1} \leq m 2d \frac{\|p\| }{\|q\| }
\leq \frac{2m d \|p\| }{(d/2)^2} = \frac{8m}{d} \|p\|  ~.
\end{equation}
Concerning $P$ we can write using also \eqref{pprime} above
\begin{equation}\label{Pprime}
\|P'\|  \leq \|p\|  + \|p'\|  \|z-b\|  \leq \|p\|  \left[ 1+
\frac{8m}{d} d \right] = (8m+1)\|p\|  ~.
\end{equation}
Consider any point $z\in K$ where $\|q\|$, and thus also $\|p\|$
is attained. We clearly have $\|P\|  \geq |P(z)|=|z-b| \|p\|$.
But here $|z-b|\geq d/5$: for in case $|z-b|\leq d/5$ we also have
$|z-a|\leq 6d/5$ by the triangle inequality, thus $|q(z)|\leq
6d^2/25 < (d/2)^2 \leq \|q\|$, as shown above. Therefore, we
conclude $\|P\|\geq (d/5) \|p\| $ and \eqref{Pprime} leads to
\begin{equation}\label{Pprimeready}
\|P'\|  \leq \frac{5(8m+1)}{d} \|P\|  < \frac{20n}{d} \|P\|
\qquad (~n:=2m+1=\deg P~)~.
\end{equation}

Now consider first the case $w>d/25$. Using $(25w/d)\geq 1$ we
obtain both for $p$ and for $P$ the estimate
\begin{equation}\label{Caseoneends}
M(p), M(P) \leq \frac{20n}{d} \leq 500 \frac{w}{d^2}n \qquad
(n:=\deg p ~\text{or} \deg P,~\text{respectively}).
\end{equation}
Note that here we have these estimates for any $n\in\NN$, without
bounds on $n$.

Let now $w< d/25$. Note that if $S$ is the strip
$S:=\{\omega=\alpha a + (1-\alpha) b + ite \in \CC~:~ 0\leq
\alpha \leq 1,~ t\in \RR\}$, then $K\subset S$, since points
outside of this strip are further than $d$ either from $a$ or
from $b$. In the following we even introduce $w^{+}:=\sup_K \Im
(\omega/e)$ and $w^{-}:=\inf_K \Im (\omega/e)$. In the current
second case of $w< d/25$, we can estimate $w^{\pm}$ by $1.02 w$.
That is, we claim that for a point $\omega=\alpha a + (1-\alpha)
b + i\beta e \in K$ with ($\alpha \in [0,1]$ and) $\beta\geq 0$,
say, we necessarily have $\beta\leq 1.02 w$. By symmetry, we may
assume that $\alpha\geq 1/2$. Put $z:=(\alpha-1/2) a +
(3/2-\alpha) b \in [a,b]$. We then find
\begin{align}\label{wwprime}
w(K)& \geq w(\{a,b,\omega\})  = w(\con \{a,b,\omega\}) \geq
w(\con\{z,\alpha a + (1-\alpha) b,\omega\}) \notag \\
& = \dist(\alpha a + (1-\alpha) b,[z,\omega]) =
\frac{(d/2)\beta}{\sqrt{(d/2)^2+\beta^2}}= \frac{\beta}{\sqrt{1 +4
(\frac{\beta}{d})^2}}~.
\end{align}
Since $\beta \leq d$ is obvious, we conclude
\begin{equation}\label{wprimerough}
\frac{d}{25}\geq w \geq \frac{\beta}{\sqrt{5}}\,,
\end{equation}
hence from \eqref{wwprime} we even have
\begin{equation}\label{wprimefine}
\beta \leq w \sqrt{1 +4 (\frac{\beta}{d})^2} \leq w
\sqrt{1+\frac{4}{125}}<1.02 w~.
\end{equation}
It follows that $w^{\pm}\leq w':=1.02 w$, as stated. Therefore,
the domain $K$ lies not only in the strip $S$, but also within the
rectangle $R:=\con\{a-iw'e,b-iw'e,b+iw'e,a+iw'e\}$. For the
central part $Q:=\{\omega\in S ~:~ |\alpha-1/2|\leq 10w/d \}$ of
$R$ we have
\begin{equation}\label{centralqprime}
\|q'\|_{K \cap Q}=\|2z-a-b\|_{K\cap Q} \leq 2 \sqrt{(10w)^2+w'^2}
< 21 w,
\end{equation}
while for the remaining part
\begin{equation}\label{sideqprime}
\|q'\|_{K\setminus Q} \leq 2d
\end{equation}
remains valid as above.

Next we estimate $q$ in $K\setminus Q$. It is easy to see that
here $\|q\|_{K\setminus Q} \leq \|q\|_{R\setminus
Q}=\left|q\left(1/2+10w/d)a+(1/2-10w/d)b+iw'e\right)\right|$,
hence
\begin{align}
\|q\|^2_{K\setminus Q} & \leq \left[\left(\frac d2+10w\right)^2
+w'^2\right] \left[\left(\frac d2 -10w\right)^2+w'^2 \right]
\notag \\
& = \left(\frac{d}{2}\right)^4 -\left(\frac{d}{2}\right)^2
\left[200w^2-2w'^2 \right] +10^4 w^4 + 200 w^2w'^2+w'^4\notag \\
& \leq \left(\frac{d}{2}\right)^4 -\left(\frac{d}{2}\right)^2
\left[197w^2 \right] +(10^4 + 210 + 2) w^4 \notag ~,
\end{align}
applying also \eqref{wprimefine}, i.e. $w'\leq 1.02 w$ in the
last step. Now $w\leq d/25$ yields
\begin{align}
\|q\|^2_{K\setminus Q} & \leq \left(\frac{d}{2}\right)^4
\left[1-788 \left(\frac{w}{d}\right)^2+ {16 \cdot
10212}\left(\frac{w}{d}\right)^4 \right]
\notag \\
& \leq \left(\frac{d}{2}\right)^4 \left[1 - 522 \left(\frac{
w}{d}\right)^2\right] \leq \left(\frac{d}{2}\right)^4 \left[1 -
\left(\frac{16 w}{d}\right)^2\right]^2 \notag ~,
\end{align}
that is, using also $\|q\|_K \geq (d/2)^2$, we find
\begin{equation}\label{qoutQ}
\frac{\|q\|_{K\setminus Q}}{\|q\|_K} \leq 1 - \left(\frac{16
w}{d}\right)^2~.
\end{equation}
Now for $z\in K\cap Q$ we have in view of \eqref{centralqprime}
\begin{align}\label{zinQ}
|p'(z)| & = m \cdot |q'(z)| \cdot |q^{m-1}(z)|  \leq 21 w m
\|q\|^{m-1} = 21wm \frac{\|p\|}{\|q\|}\notag \\ & \leq \frac{21wm
\|p\|}{(d/2)^2}=42\frac{w}{d^2}n\|p\|~,
\end{align}
and for $z\in K\setminus Q$ using $\|q\|_K \geq (d/2)^2$,
\eqref{sideqprime} and \eqref{qoutQ} we get
\begin{align}\label{znotinQ}
|p'(z)| & \leq m\cdot 2d \cdot \| q\|^{m-1}_{K\setminus Q} \leq
2md \frac{\|p\|}{\|q\|} \left[1 - \left(\frac{16
w}{d}\right)^2\right]^{m-1} \notag \\ & \leq \frac{8m}{d} \|p\|
\left[1 - \left(\frac{16 w}{d}\right)^2\right]^{m-1}~.
\end{align}
Now in view of $w<d/25$, a standard calculation shows that
\begin{equation}\label{expoest}
\left[1 - \left(\frac{16 w}{d}\right)^2\right]^{m-1} \leq \frac{25
w}{d} \qquad \text{if}~~~ m \geq m_0:=\left( \frac{d}{16w}
\right)^2 \log \left( \frac{d}{16w} \right)~.
\end{equation}
Indeed, as $\log(1-x) < - x$ for all $0<x<1$, using $w<d/25$ we
find
$$
(m-1) \log \left[1-\left(\frac{16w}{d}\right)^2\right] < -(m-1)
\left(\frac{16w}{d}\right)^2 < -m \left(\frac{16w}{d}\right)^2 +
0.41,
$$
which entails for $m \geq m_0$ that
$$
\left[1 - \left(\frac{16 w}{d}\right)^2\right]^{m-1} < e^{ - m_0
\left(\frac{16w}{d}\right)^2 + 0.41} =e^{-\log
\left(\frac{d}{16w}\right) + 0.41} < \frac{25w}{d}~.
$$
It follows from \eqref{znotinQ} and \eqref{expoest} that
\begin{equation}\label{pprimeoutQ}
\|p'\|_{K\setminus Q} \leq 100 \frac{w}{d^2} n \|p\|~.
\end{equation}
Collecting \eqref{zinQ} and \eqref{pprimeoutQ} we get also in
this case of $w<d/25$ the estimate
\begin{equation}\label{psmallw}
\|p'\| \leq 100 \frac{w}{d^2} n \|p\|\qquad (n=2m=\deg p,~~ m\geq
m_0)\, .
\end{equation}
It remains to consider the odd degree case of $n=2m+1$, i.e. $P$.
Now write
\begin{equation}\label{Pprimez}
|P'(z)|\leq |p(z)|+|p'(z)|\cdot|z-b| \leq |p(z)|+d \|p'\| \leq
(1+100 \frac{w}{d} 2m ) \|p\|\quad (m\geq m_0),
\end{equation}
in view of \eqref{psmallw}. As shown above, we have $\|P\|\geq
\|p\|/(d/5)$, while $m\geq m_0$ entails $1\leq m/m_0 < m(w/d)
(16^2/25( (1/\log(25/16)) < 30 mw/d $, hence \eqref{Pprimez}
yields
$$
\|P'\| \leq \frac{230mw}{d} \|p\| \leq \frac{1150mw}{d^2} \|P\|~.
$$
Since now $n=2m+1>2m$, we finally find
\begin{equation}\label{Pprimenorm}
\|P'\| < 600 \frac{w}{d^2} n \|P\|\qquad (n=2m+1=\deg
P,~~m>m_0)\,.
\end{equation}
\end{proof}

%%% \newpage

\mysection{Acknowledgements and comments}

\bigskip

Because \cite{Rev} will not be published in a journal, a full,
self-contained proof was presented here. At the same time, this
was meant to provide also a clear explanation and documentation
of the origin and development of the various ideas that have led
to the result.

The author is indebted to G\'abor Hal\'asz for his generous
contribution of an essential idea, explained at the beginning of
\S \ref{sec:proof} above. Although he did not accept being
included as coauthor, the paper contains, in fact, a joint result
with him.

%%%%%%%%%%%%%%%%%%%%%%%%%%%%%%%%%%%%%%%%%%%%%
% REFERENCES
%%%%%%%%%%%%%%%%%%%%%%%%%%%%%%%%%%%%%%%%%%%%%
\newpage

\noindent
\ \\
{\bf Bibliography}
\\

\bigskip

\noindent
{\sc\small Alfr\' ed R\'enyi Institute of Mathematics, \\
Hungarian Academy of Sciences, \\
Budapest, POB 127, \\ 1364 Hungary}\\
E-mail: {\tt revesz@renyi.hu}

\end{document}